\documentstyle{article}

\title{Harmonic maps between generalized\\
Lagrange spaces}
\author{Mircea Neagu}
\date{}
\begin{document}
\maketitle
\begin{abstract}
In Section 1 the author defines the notion of harmonic map between generalized
Lagrange spaces. Section 2 analyses the particular case when the
generalized Lagrange spaces are Lagrange spaces of electrodynamics. In Section 3
it is proved that for certain systems of differential or partial differential
equations, the solutions are harmonic maps between certain generalized Lagrange spaces,
in the sense of Section 1.
Section 4 describes the main properties of the generalized Lagrange spaces
constructed in Section 3.
\end{abstract}
\section{Introduction}

\hspace{5mm}Let $(M^m,\;g_{\alpha\beta})$ and $(N^n,\;h_{ij})$ be two generalized
Lagrange spaces, where $m$, respectively $n$, is the dimension of $M$,
respectively $N$. The manifold $M$, respectively $N$, is coordinated by
$(a^\alpha)$, respectively $(x^i)$.
On $M\times N$, the first $m$ coordinates are indexed by
$\alpha,\beta,\gamma,\ldots$ and the last $n$ coordinates are indexed by
$i,j,k,\ldots$.
The fundamental metric tensors are expressed locally by

i) $g_{\alpha\beta}=g_{\alpha\beta}(a,b),\;\forall\;\alpha ,\beta=
\overline {1,m}$, where $(a,b)=(a^\mu ,b^\mu )$ are adapted coordinates
on $TM$.

ii) $h_{ij}=h_{ij}(x,y),\;\forall\;i,j=\overline{1,n}$, where
$(x,y)=(x^k,y^k)$ are adapted coordinates on $TN$.

On $M\times N$, we consider an arbitrary tensor of type $(1,2)$, denoted by $P$, with all components
null except $P^\beta_{\alpha i}(a, x)$ and $P^j_{\alpha i}(a, x)$, where $\alpha ,\beta =
\overline{1,m},\;i,j=\overline{1,n}$, which will be called
{\it tensor of connection}. The connection tensor $P$ allows one to build the
directions $b$ and $y$ of the metric tensors $g_{\alpha\beta}(a,b)$ and
$h_{ij}(x,y)$ used in the construction of the energy functional $E$ whose extremals
will be the harmonic maps between the generalized Lagrange spaces $M$ and $N$.
In conclusion, this tensor makes the connection between the metric structures
of the spaces $M$, respectively $N$, and the harmonic maps that we will define.

We assume that the manifold $M$ is connected, compact,
orientable and endowed with a Riemannian metric $\varphi_{\alpha\beta}$.
These conditions assure the existence of a volume element and, implicitly,
of a theory of integration on $M$. Using the generalized Lagrange metrics
$g_{\alpha\beta}$ and $h_{ij}$ and
the tensor of connection $P$, we can define the following
$\left(
\begin{array}{ccc}
&P&\\
g&\varphi&h
\end{array}\right)$
{\it-energy functional},
$$
E=E^{\scriptstyle P}_{g\varphi h}:C^\infty (M,N)\to R,
$$
$$
E^{\scriptstyle P}_{g\varphi h}(f)=\displaystyle{{1\over 2}\int_M
g^{\alpha\beta}(a,b(a,f^k,f^k_\gamma))
h_{ij}(f(a),y(a,f^k,f^k_\gamma))f^i_\alpha f^j_\beta\sqrt\varphi\;da,}
$$
where
$\left\{\begin{array}{l}\medskip
\displaystyle{
f^i=x^i(f),\;f^i_\alpha ={\partial f^i\over\partial a^\alpha},\;
\varphi=det(\varphi_{\alpha\beta})},\\\medskip
\displaystyle{b(a,f^k,f^k_\gamma)=b^\gamma (a)
\left.{\partial\over\partial a^\gamma}\right\vert_a
\stackrel{\hbox{def}}=\varphi^{\alpha\beta}(a)f^i_\alpha (a)P^\gamma_{\beta i}(a,f(a))
\left.{\partial\over\partial a^\gamma}\right\vert_a
},\\
\displaystyle{
y(a,f^k,f^k_\gamma)=y^k(a)\left.{\partial\over\partial x^k}\right\vert_{f(a)}
\stackrel{\hbox{def}}=\varphi^{\alpha\beta}(a)f^i_\alpha (a)P^k_{\beta i}(a,f(a))
\left.{\partial\over\partial x^k}\right\vert_{f(a)}}
.\end{array}\right.$

{\bf Definition.}
A map $f\in C^\infty(M,N)$ is
$\left(
\begin{array}{ccc}
&P&\\
g&\varphi&h
\end{array}\right)$
{\it-harmonic} iff $f$ is a critical point for the functional
$E^{\scriptstyle P}_{g\varphi h}$.

{\bf Particular cases.} i) If $g_{\alpha\beta}(a,b)=\varphi_{\alpha\beta}(a)$
and $h_{ij}(x,y)=h_{ij}(x)$
are Riemannian me\-trics and the connection tensor is an arbitrary one, it
recovers the classical definition of a harmonic map
between two Riemannian manifolds \cite{1}. We remark that, in this case, the
definition of harmonic maps is independent of the connection tensor field $P$.

ii) If we take $N=R,\;h_{11}=1$ and the tensor of connection is of the form
$P=(\delta^\alpha_\beta,P^1_{\beta 1})$,
we obtain $C^\infty(M,N)={\cal F}(M)$ and the energy functional becomes
$$
E^{\scriptstyle P}_{g\varphi 1}(f)=\displaystyle{
{1\over 2}\int_Mg^{\alpha\beta}(a,grad_\varphi f)f_\alpha f_\beta\sqrt\varphi
\;da\;,\;\forall f\in{\cal F}(M).}
$$

iii) If we consider $M=[a,b]\subset R,\;\varphi_{11}=g_{11}=1$ and the
tensor of connection is\linebreak
$P=(P^1_{1i},\delta^k_i)$,
we obtain $C^\infty (M,N)=\{x:[a,b]\to N\;\vert\;x-C^\infty\;\hbox{differentiable}
\}$. Denoting $C^\infty(M,N)
=\Omega_{a,b}(N)$, the energy functional should be
$$
E^{\scriptstyle P}_{11h}(x)=\displaystyle{{1\over 2}\int_a^bh_{ij}(x(t),\dot x(t)){dx^i\over
dt}{dx^j\over dt}dt,\;\forall x\in\Omega_{a,b}(N).}
$$
In conclusion, the
$\left(
\begin{array}{ccc}
&P&\\
1&1&h
\end{array}\right)$-harmonic curves are exactly the geodesics of the generalized
Lagrange space $(N,h_{ij}(x,y))$ \cite {2}.

\section{Harmonic maps between Lagrange spaces of electrodynamics}

\hspace{5mm}Let $(M^m,L_M)$ and $(N^n,L_N)$ be Lagrange spaces with
the Lagrangians
$$
L_M(a,b)=g_{\alpha\beta}(a)b^\alpha b^\beta+g_{\alpha\beta}(a)U^\alpha (a)b^
\beta +F(a),
$$
$$L_N(x,y)=h_{ij}(x)y^iy^j+h_{ij}(x)V^i(x)y^j+G(x),
$$
where
\begin{itemize}
\item $g_{\alpha\beta}$ (resp. $h_{ij}$) is Riemannian metric on $M$ (resp. $N$)
representing the gravitational potentials on $M$ (resp. $N$).
\item $U^\alpha$ (resp. $V^i$) is a vector field on $M$ (resp. $N$) representing
the electromagnetic potentials.
\item $F$ (resp. $G$) is a smooth function on $M$ (resp. $N$) representing the
potential function.
\end{itemize}

The fundamental metric tensors of these Lagrangians are
$$
\displaystyle{g_{\alpha\beta}(a)=
{\partial^2L_M\over\partial b^\alpha\partial b^\beta},\;
h_{ij}(x)={\partial^2L_N\over\partial y^i\partial y^j}.}
$$
Taking an arbitrary tensor of connection $P$, the energy functional becomes
$$
E_{g\varphi h}(f)=\displaystyle{
{1\over 2}\int_Mg^{\alpha\beta}(a)h_{ij}(f(a))f^i_\alpha f^j_\beta\sqrt\varphi\;da}.
$$
We remark that, in this case, the energy functional is independent of the
connection tensor $P$.

The Euler-Lagrange equations will be, obviously, the equations of harmonic maps,
namely
$$\displaystyle{
g^{\alpha\beta}\left\{f^k_{\alpha\beta}-\left [G^\gamma_{\alpha\beta}+{1\over 2}
{\partial\over\partial a^\alpha}\left(\ln{g\over\varphi}\right)\delta^\gamma_
\beta\right]f^k_\gamma+H^k_{ij}f^i_\alpha f^j_\beta\right\}=0\;,\;\forall\;k=
\overline{1,n},}
$$ where
\begin{itemize}
\item $f^k_{\alpha\beta}=\displaystyle{{\partial^2f^k\over\partial a^\alpha\partial
a^\beta}\;,\;g=det(g_{\alpha\beta})\;,\;\varphi =det(\varphi_{\alpha\beta})}$.
\item $G^\gamma_{\alpha\beta}$ are the Christoffel symbols of the metric
$g_{\alpha\beta}$.
\item $H^k_{ij}$ are the Christoffel symbols of the metric $h_{ij}$.
\end{itemize}

{\bf Remarks.} i) If $g_{\alpha\beta}=\varphi_{\alpha\beta}$, we recover the
classical equations of harmonic maps between two Riemannian manifolds.

ii) Denoting $\displaystyle{\Delta^\gamma_{\alpha\beta}=
G^\gamma_{\alpha\beta}+
{1\over 2}{\partial\over\partial a^\alpha}\left(\ln{g\over\varphi}\right)\delta^\gamma_
\beta}$, we remark that $\Delta^\gamma_{\alpha\beta}$ represent the components
of a linear connection induced by the metrics $g_{\alpha\beta}$ and $\varphi_
{\alpha\beta}$.

By the last remark, we can give the following

{\bf Definition.} Let $g,\varphi$ be Riemannian metrics on $M$. The curve
$c:I\to M$, expressed locally by $c(t)=(a^\alpha (t))$, is a $(g,\varphi)$
-{\it geodesic} iff $c$ is an autoparallel curve of the connection $\Delta^\gamma_
{\alpha\beta}$ induced by the metrics $g$ and $\varphi$, namely
$$\displaystyle{{d^2a^\gamma\over dt^2}=-\Delta^\gamma_{\alpha\beta}
{da^\alpha\over dt}{da^\beta\over dt}}.
$$

{\bf Remarks.} i) If $g=\varphi$, then we recover the classical definition of
a geodesic on the Riemannian manifold $(M,g=\varphi)$.

ii) It is obviously that a $(g,\varphi)$-geodesic is a reparametrized
geodesic of the metric\nolinebreak $\;g$.

{\bf Theorem.} {\it
Let $f:(M,L_M)\to (N,L_N)$ be a smooth map which carries
$(g,\varphi)$-geodesics into $h$-geodesics. Then $f$ is harmonic map.}

{\bf Proof.} Let $c:I\subset R\to M\;,\;c(t)=(a^\alpha(t))$ be a
$(g,\varphi)$-geodesic. Then we have
$\displaystyle{{d^2a^\gamma\over dt^2}=-\Delta^\gamma_{\alpha\beta}{da^\alpha
\over dt}{da^\beta\over dt}}$. Because $\overline{c}(t)=f(c(t))$ is $h$-geodesic,
it follows that $\displaystyle{{d^2\overline{c}^k\over dt^2}+H^k_{ij}
{d\overline{c}^i\over dt}{d\overline{c}^j\over dt}=0}$. But $\displaystyle{
{d\overline{c}^k\over dt}=f^k_\alpha(c(t)){da^\alpha\over dt}\Rightarrow
{d^2\overline{c}^k\over dt^2}=f^k_{\alpha\beta}{da^\alpha\over dt}{da^\beta
\over dt}+f^k_\alpha{d^2a^\alpha\over dt^2}}$.
In conclusion we obtain
$$\displaystyle{{d^2a^\alpha\over dt^2}f^k_\alpha+f^k_{\alpha\beta}{da^\alpha\over dt}
{da^\beta\over dt}+H^k_{ij}f^i_\alpha f^j_\beta{da^\alpha\over dt}{da^\beta\over
dt}=0\Rightarrow}
$$
$$\displaystyle{\left(f^k_{\alpha\beta}-\Delta^\gamma_{\alpha\beta}f^k_\gamma+
H^k_{ij}f^i_\alpha f^j_\beta\right){da^\alpha\over dt}{da^\beta\over dt}=0\;,
\;\forall k=\overline{1,n}}\;\Box.
$$

\section{Geometrical interpretation of solutions of certain PDEs
of order one}

\hspace{5mm}The problem of finding a geometrical structure of Riemannian type on a manifold
$M$ such that the orbits of an arbitrary vector field $X$ should be geodesics, was
intensively studied by Sasaki. The results were not satisfactory, but,
in his study, Sasaki
discovered the well known almost contact structures on a manifold of odd\linebreak
dimension \cite{6}. After the introduction of generalized Lagrange spaces
by Miron \cite{2}, the same problem is resumed by
Udri\c ste \cite{8,9}. This succeded to discover a Lagrange structure on $M$,
depending of the vector field $X$ and an associated (1,1)-tensor field,
such that the orbits of $C^2$ class should be geodesics.
Moreover, he formulated a more general problem \cite{9}, namely

1) Are there structures of Lagrange type such that the solutions of
certain PDEs of order one should be {\it harmonic maps}?

2) What is a {\it harmonic map} between two generalized Lagrange spaces?

A partial answer of these questions is offered by author in his paper \cite{4},
using the notion of harmonic map on a direction between a Riemannian manifold
and a generalized Lagrange manifold. The notion of {\it harmonic map} between
two generalized Lagrange spaces introduced in \cite{7}
allows one to extend the results of previous papers \cite{4}, \cite{8},
\cite{9} and
to obtain a beautiful geometrical interpretation of the solutions of the certain
PDEs of order one.

For every smooth map $f\in C^\infty(M,N)$, we use the following notation
$$\displaystyle{
\left.\delta f=f^i_\alpha da^\alpha\vert_a
\otimes{\partial\over\partial y^i}\right\vert_{f(a)}\in \Gamma(T^*M\otimes
f^{-1}(TN))}.$$
On $M\times N$, let $T$ be one tensor of type $(1,1)$ with all
components equal to zero except $(T^i_{\alpha} )_{i=\overline{1,n}\\
\atop
\alpha=\overline{1,m}}$.
\medskip
Let the system of partial
differential equations
$$\displaystyle{
\delta f=T\;\hbox{expressed locally by}\;{\partial f^i\over\partial a^\alpha}
=T^i_\alpha(a,f(a)).
}\leqno(E)$$

If $(M,\varphi_{\alpha\beta})$ and $(N,\psi_{ij})$ are Riemannian
manifolds, we can build a scalar product on $\Gamma(T^*M\otimes f^{-1}(TN))$ by
$<T,S>=\varphi^{\alpha\beta}(a)\psi_{ij}(f(a))T^i_\alpha S^j_\beta$, where
$\displaystyle{T=T^i_\alpha da^\alpha\otimes{\partial\over\partial y^i}}$ and
$\displaystyle{S=S^j_\beta da^\beta\otimes{\partial\over\partial y^j}}$.

Under these conditions, we can prove the following

{\bf Theorem.} {\it
If $(M,\varphi),(N,\psi)$ are Riemannian manifolds and $f\in C^\infty(M,N)$ is a
solution of the system $(E)$, then $f$ is a solution of the variational problem
asociated to the functional
${\cal L}_T:C^\infty(M,N)\backslash\{f\;\vert\;\exists\;a\in M\;\hbox{such that}\;
<\delta f,T>(a)=0\}\to R_+$,
$$
\displaystyle{{\cal L}_T(f)={1\over 2}\int_M{\Vert\delta f\Vert^2\Vert T\Vert^2
\over <\delta f,T>^2}\sqrt{\varphi}\;da={1\over 2}\int_M{\Vert T\Vert^2\over <
\delta f, T>^2}\varphi^{\alpha\beta}\psi_{ij}f^i_\alpha f^j_\beta\sqrt{\varphi}\;da}.
$$}

{\bf Proof.} In the space $\Gamma(T^*M\times
f^{-1}(TN))$, the Cauchy inequality for the scalar product $<,>$ holds.
It follows that the following inequality is true,
\linebreak$<T,S>^2\leq\nolinebreak\Vert T\Vert^2\Vert S\Vert^2\;,\;\forall\;
T,S\in\Gamma(T^*M\times\nolinebreak f^{-1}(TN))$, with equality if
and only if there exists
${\cal K}\in{\cal F}(M)$ such that $T={\cal K}S$. Consequently,
for every $f\in C^\infty (M,N)$ we have
$$\displaystyle{{\cal L}(f)={1\over 2}\int_M{\Vert\delta f\Vert^2
\Vert T\Vert^2\over<\delta f,T>^2}\sqrt{\varphi}\;da\geq{1\over 2}\int_M\sqrt
{\varphi}\;da={1\over 2}Vol_\varphi(M)}.
$$
Obviously, if $f$ is a solution of the system $(E)$, we obtain
$\displaystyle{{\cal L}_T(f)={1\over 2}Vol_\varphi(M)}$, that is, $f$ is a
global minimum point for ${\cal L}_T$ $\Box$.

{\bf Remarks.} i) In certain particular cases of the system $(E)$, the
functional ${\cal L}_T$ becomes
exactly a functional of type
$\left(
\begin{array}{ccc}
&P&\\
g&\varphi&h
\end{array}\right)$-energy.

ii) The global minimum points of the functional ${\cal L}_T$ are solutions of
the system $\delta f={\cal K}T$, where ${\cal K}\in{\cal F}(M)$, not necessarily
with ${\cal K}=1$.

iii) Replacing the Riemannian metric $\psi_{ij}$ by a pseudo-Riemannian
metric, the preceding theorem survives because the form of the Euler-Lagrange
equations remains unchanged. The difference is that the
solutions of the system $(E)$, in the pseudo-Riemannian case, are not the
global minimum points for the functional ${\cal L}_T$. Moreover, the statement (ii),
of above, does not hold.

{\bf Fundamental examples.}

{\bf 1. Orbits}

For $M=([a,b],1)$ and $T=\xi\in\Gamma(x^{-1}(TN))$, the system
$(E)$ becomes
$$\displaystyle{{dx^i\over dt}=\xi^i(x(t)),\;x:[a,b]\to N},\leqno{(E_1)}$$
that is the system of orbits for $\xi$, and
the functional ${\cal L}_\xi$ is
$$\displaystyle{{\cal L}_\xi(x)={1\over 2}\int^b_a{\Vert\xi\Vert^2_\psi\over
[\xi^b(\dot x)]^2}\psi_{ij}{dx^i\over dt}{dx^j\over dt}dt},
$$
where $\xi^b=\xi_idx^i=\psi_{ij}\xi^jdx^i$.
Hence the functional ${\cal L}_\xi$ is a
$\left(
\begin{array}{ccc}
&P&\\
1&1&h
\end{array}\right)$-energy (see (iii) of first particular cases of this paper),
where
$$h_{ij}:TN\backslash\{y\;\vert\;\xi^b(y)=0\;\hbox{for some}\;y\}\to R$$
is defined by
$$\displaystyle{
h_{ij}(x,y)={\Vert\xi\Vert^2_\psi\over[\xi^b(y)]^2}\psi_{ij}(x)=\psi_{ij}(x)
\exp{\displaystyle{\left[
2\ln{\Vert\xi\Vert_\psi\over\vert\xi^b(y)\vert}\right]} }
}.$$
This case is studied in other way by Udri\c ste in \cite{8}-\cite{9}.

{\bf 2. Pfaffian systems}

For $N=(R,1)$ and $T=A\in\Lambda^1(T^*M)$, the system $(E)$
becomes
$$df=A,\;f\in{\cal F}(M),\leqno{(E_2)}$$
that is a Pfaffian system, and the functional ${\cal L}_T$ reduces to
$$
\displaystyle{{\cal L}_A(f)={1\over 2}\int_M{\Vert A\Vert^2_\varphi\over[A(
grad_\varphi f)]^2}\varphi^{\alpha\beta}f_\alpha f_\beta\sqrt{\varphi}da}.
$$
Hence, the functional ${\cal L}_A$ is a
$\left(
\begin{array}{ccc}
&P&\\
g&\varphi&1
\end{array}\right)$-energy
(see (ii) of first particular cases of this paper), where
$g_{\alpha\beta}:TM\backslash\{b\;\vert\; A(b)=0\;\hbox{for some}\;b\}\to R$ is defined by
$$\displaystyle{
g_{\alpha\beta}(a,b)={[A(b)]^2\over\Vert A\Vert^2_\varphi}
\varphi_{\alpha\beta}(a)
=\varphi_{\alpha\beta}(a)\exp{
\displaystyle{
\left[2\ln{\vert A(b)\vert\over\Vert A\Vert_\varphi}
\right]}}
}.$$

{\bf 3. Pseudolinear functions}

We suppose that $T^k_\beta (a,x)=\xi^{k}(x)A_\beta (a)$,
where $\xi^k$ is vector field on $N$ and $A_\beta$ is 1-form on $M$.
In this case the system $(E)$ is
$$\displaystyle{{\partial f^k\over\partial a^\beta}=\xi^k(f)A_\beta (a)}
\leqno{(E_3)}
$$ and the functional ${\cal L}_T$ is expressed by
$$\displaystyle{{\cal L}_T(f)={1\over 2}\int_M{\Vert\xi\Vert^2_\psi\Vert A
\Vert^2_\varphi\over[A(b)]^2}\varphi^{\alpha\beta}\psi_{ij}f^i_\alpha f^j_\beta
\sqrt\varphi da=
}$$
$$\displaystyle{={1\over 2}\int_Mg^{\alpha\beta}(a,b)h_{ij}(f(a))f^i_\alpha f^j
_\beta\sqrt\varphi da},
$$
where $h_{ij}(x)=\Vert\xi\Vert^2_\psi\psi_{ij}(x)$, the tensor of connection is
$P^\gamma_{i\beta}(x)=\delta^\gamma_\beta\xi_i(x),\linebreak b^\gamma=
\varphi^{\alpha\beta}f^i_\alpha P^\gamma_{i\beta}$
and the Lagrange metric tensor
$$g_{\alpha\beta}:TM\backslash\{b\;\vert\;A(b)=0\;\hbox{for some}\;b\}\to R$$
is defined by
$$\displaystyle{g_{\alpha\beta}(a,b)={[A(b)]^2\over\Vert A\Vert^2_\varphi}
\varphi_{\alpha\beta}(a)=\varphi_{\alpha\beta}(a)
\exp{\displaystyle{\left[2\ln{\vert A(b)\vert\over\Vert A\Vert_\varphi}\right]}}
}.$$
It follows
that the functional ${\cal L}_T$ becomes a
$\left(
\begin{array}{ccc}
&P&\\
g&\varphi&h
\end{array}\right)$-energy.\medskip

{\bf Remark.} Take $M$ to be an open subset in $(R^n,\varphi=\delta)$ and
$N=(R,\psi=1)$, the system from the third example is
$$\displaystyle{{\partial f\over\partial a^\alpha}=\xi(a)A_\alpha(f(a)),\;
\forall\;\alpha=\overline{1,m}.}\leqno{(PL)}
$$
Supposing that $(grad\;f)(a)\ne 0,\;\forall\;a\in M$, the solutions of this
system are the well known {\it pseudolinear functions} \cite{5}. These functions
have the following property,

$-$for every fixed point $x_0\in M$, the hypersurface of constant level
$$M_{f(x_0)}=\{x\in M\;\vert\;f(x)=f(x_0)\}$$ is totally geodesic
\cite {5} (i. e. the second
fundamental form vanishes identically).\medskip

In conclusion, the pseudolinear functions are examples of harmonic maps
between the generalized Lagrange spaces
$\displaystyle{\left(M,g_{\alpha\beta}(a,b)={[A(b)]^2\over\Vert A\Vert^2}
\delta_{\alpha\beta}\right)}$ and
$\displaystyle{(R,h(x)=\xi^2(x))}$.
For example, we have the following pseudolinear
functions \cite{5}:\medskip

{\bf 3. 1.} $f(a)=e^{<v,a>+w}$, where $v\in M,\;w\in R$, is
solution for the system $(PL)$ with $\xi(a)=1$ and
$A(f(a))=f(a)v$.\medskip

{\bf 3. 2.} $\displaystyle{f(a)={<v,a>+w\over<v',a'>+w'}}$, where $v,v'\in M
\;,\;w,w'\in R$, is solution for $(PL)$ with $\displaystyle{\xi(a)={1\over<v',a>+w'}}$
and $A(f(a))=v-f(a)w$.

{\bf Remark.} The preceding cases appear also in \cite{7} and, from another point
of view, in \cite{4}. The following case is the main novelty of this paper.

{\bf 4. The general case}

If we have $T^i_\alpha(a,x)=\sum_{r=1}^t\xi_r^i(x)A^r_\alpha(a)$, where
$\{\xi_r\}_{r=\overline{1,t}}\subset{\cal X}(N)$ is a family of vector
fields on $N$ and $\{A^r\}_{r=\overline{1,t}}\subset\Lambda^1(T^*M)$ is a
family of 1-forms on $M$, the system of equations $(E)$ reduces to
$$\displaystyle{{\partial f^i\over\partial a^\alpha}=\sum_{r=1}^t\xi_r^i(f)A^r_\alpha(a)
}.\leqno{(E_4)}$$

Without loss of generality, we can suppose that
$\{\xi_r\}_{r=\overline{1,t}}\subset{\cal X}(N)$
(resp. $\{A^r\}_{r=\overline{1,t}}\subset\Lambda^1(T^*M)$) are linearly independent.
In these conditions, we shall have $t\leq\min\{m,n\}$, where $m=dim\;M$ and
$n=dim\;N$.

{\bf 4. 1.} Assume that $\{\xi_r\}_{r=\overline{1,t}}\subset{\cal X}(N)$ is
an orthonormal system of vector fields with respect to the Riemannian metric
$\psi_{ij}$ on $N$. Let $B\in\Lambda^1(T^*M)$ be an arbitrary unit 1-form
on $M$. With our assumptions, by a simple calculation, we obtain
$$\Vert T\Vert^2=\varphi^{\alpha\beta}\psi_{ij}\xi^i_rA^r_\alpha\xi^j_sA^s_\beta=
\sum^t_{r,s=1}<\xi_r,\xi_s>_\psi<A^r,A^s>_\varphi=\sum_{r=1}^t\Vert A^r\Vert^2_\varphi,$$
$$<\delta f, T>=\varphi^{\alpha\beta}\psi_{ij}f^i_\alpha\xi^j_rA^r_\beta
B^\mu B_\mu.$$
Defining the tensor of connection by
$P^\gamma_{i\beta}(a,x)=\psi_{ij}(x)\xi^j_r(x)A^r_\beta(a)B^\gamma(a)$
and \linebreak $b^\gamma=\varphi^{\alpha\beta}f^i_\alpha P^\gamma_{i\beta}$, the
functional ${\cal L}_T$ takes the form
$$\displaystyle{{\cal L}_T(f)={1\over 2}\int_M{\sum_{r=1}^t\Vert A^r\Vert^2_\varphi\over[B(b)]^2}
\varphi^{\alpha\beta}\psi_{ij}f^i_\alpha f^j_\beta\sqrt{\varphi}\;da=
{1\over 2}\int_Mg^{\alpha\beta}(a,b)\psi_{ij}(f(a))f^i_\alpha f^j_\beta\sqrt
{\varphi}\;da
},$$
where the Lagrange metric tensor
$g_{\alpha\beta}:TM\backslash\{b\;\vert\;B(b)=0\;\hbox{for some}\;b\} \to R$ is expressed by
$$\displaystyle{
g_{\alpha\beta}(a,b)={[B(b)]^2\over\sum_{r=1}^t\Vert A^r\Vert^2_\varphi}
\varphi_{\alpha\beta}(a)=\varphi_{\alpha\beta}(a)
\exp{\displaystyle{\left[2\ln{\vert B(b)\vert\over
\sqrt{\sum_{r=1}^t\Vert A^r\Vert^2_\varphi}}\right]}}
}.$$
Consequently, the functional ${\cal L}_T$ is a
$\left(
\begin{array}{ccc}
&P&\\
g&\varphi&\psi
\end{array}\right)$-energy.

{\bf 4. 2.} As above, we assume that the system $\{A^r\}_{r=\overline{1,t}}$
of 1-forms is orthonormal with respect to the metric $\varphi^{\alpha\beta}$
and we choose an arbitrary unit vector field $X\in{\cal X}(N)$. By analogy to
{\bf 4. 1} we shall have
$$\Vert T\Vert^2=\sum_{r=1}^t\Vert\xi_r\Vert^2_\psi\;\hbox{and}\;
<\delta f,T>=\varphi^{\alpha\beta}\psi_{ij}f^i_\alpha\xi^j_rA^r_\beta X^kX_k.
$$
Using the notations
$P^k_{i\beta}(a,x)=\psi_{ij}(x)\xi^j_r(x)A^r_\beta(a)X^k(x)$ and
$y^k=\varphi^{\alpha\beta}f^i_\alpha P^k_{i\beta}$ we \linebreak obtain
the following expression of the functional ${\cal L}_T$,
$$\displaystyle{{\cal L}_T(f)={1\over 2}\int_M{\sum_{r=1}^t\Vert\xi_r\Vert^2
_\psi\over[X^b(y)]^2}\varphi^{\alpha\beta}\psi_{ij}f^i_\alpha f^j_\beta\sqrt
{\varphi}\;da=
{1\over 2}\int_M\varphi^{\alpha\beta}(a)h_{ij}(f(a),y)f^i_\alpha f^j_\beta
\sqrt{\varphi}\;da
},$$
where the Lagrange metric tensor
$h_{ij}:TN\backslash\{y\;\vert\;X^b(y)=0\;\hbox{for some}\;y\} \to R$ is
$$\displaystyle{
h_{ij}(x,y)={\sum_{r=1}^t\Vert\xi_r\Vert^2_\psi\over[X^b(y)]^2}
\psi_{ij}(x)=\psi_{ij}(x)\exp{\displaystyle{
\left[2\ln{\sqrt{\sum_{r=1}^t\Vert\xi_r\Vert^2_\psi}\over\vert
X^b(y)\vert}\right]}}
}.$$
Obviously, the functional ${\cal L}_T$ is a
$\left(
\begin{array}{ccc}
&P&\\
\varphi&\varphi&h
\end{array}\right)$-energy.

{\bf Remarks.} i) For the use of the above, we assume {\it a priori} a Riemannian
metric $\varphi$ or $\psi$ on $M$ or $N$ such that the system of covectors
$\{A^r\}_{r=\overline{1,t}}$ or of vectors $\{\xi_r\}_{r=\overline{1,t}}$
is orthonormal. This fact is always possible. In conclusion, our assumptions
on the orthonormality of these systems do not restrict the generality of problem.

ii) In the particular case, $n=r=1\;,\;N=R\;,\;\psi_{11}=1\;,\;\xi=
\displaystyle{d\over dx}$, and $A\in\Lambda^1(T^*M)$ is an arbitrary 1-form
on $M$,
we recover the Pfaffian system $df=A$. In this situation, taking
$B\in\Lambda^1(T^*M)$ to be an arbitrary unit 1-form, we obtain, for the functional
${\cal L}_T$, the expression
$$\displaystyle{{\cal L}_T(f)={1\over 2}\int_M{\Vert A\Vert^2\over[B(b)]^2}
\varphi^{\alpha\beta}f_\alpha f_\beta\sqrt{\varphi}\;da={1\over 2}\int_M
{\Vert A\Vert^2\over[A(grad_\varphi f)]^2}\varphi^{\alpha\beta}f_\alpha f_\beta
\sqrt{\varphi}\;da
},$$ where the tensor of connection is defined by $P_{1\beta}^\gamma=A_\beta
B^\gamma$ and the direction $b$ is\linebreak
$b^\gamma=\varphi^{\alpha\beta}f^i_\alpha A_\beta B^\gamma=B^\gamma
A(grad_\varphi f)$.

Since the last integral is the functional ${\cal L}_A$ from the example {\bf 1},
we remark that the solutions of the Pfaffian system $df=A$ can be regarded in an
infinity manner as
$\left(
\begin{array}{ccc}
&P&\\
g&\varphi&1
\end{array}\right)$-harmonic maps. This fact appears because the tensor of
connection $P$ and the generalized Lagrange metric $g$ are dependent of
the arbitrary unit covector field $B$.

iii) Analog to ii), if we take $m=r=1\;,\;M=[a,b]\;,\;
g_{11}=\varphi_{11}=1\;,\; A=dt$ and $\xi\in{\cal X}(N)$ is an arbitrary
vector field on $N$, we find the system of orbits for $\xi$, that is,
$$\displaystyle{{dx^i\over dt}=\xi^i(x(t)),\;x:[a,b]\to N
.}$$
Starting with $X\in{\cal X}(N)$ an arbitrary unit vector field, the functional
${\cal L}_T$ becomes
$$\displaystyle{{\cal L}_T(x)={1\over 2}\int_a^b{\Vert\xi\Vert^2_\psi\over[X^b(y)]^2}
\psi_{ij}{dx^i\over dt}{dx^j\over dt}dt={1\over 2}\int_a^b{\Vert\xi\Vert^2_\psi
\over[\xi^b(\dot x)]^2}\psi_{ij}{dx^i\over dt}{dx^j\over dt}dt
},$$
where the tensor of connection is $P_{i1}^k=\xi_i^bX^k$ and $y^k=\xi^b(\dot x)X^k$.
Because we can vary the tensor of connection $P$ and the generalized Lagrange
metric $h$ by the arbitrary unit vector field $X$, we remark that the trajectories of
the vector field $\xi$ can be also regarded in an infinity manner as
$\left(
\begin{array}{ccc}
&P&\\
1&1&h
\end{array}\right)$-harmonic maps.

\section{Lagrange geometry asociated to PDEs of order one}

\hspace{5mm}We first remark that, in all above cases, the solutions
of $C^2$ class of the system $\delta f=T$
becomes harmonic maps between generalized Lagrange spaces, in the sense defined
in this paper. Moreover, the above generalized Lagrange structures are of type
$(M^n,e^{2{\sigma(x,y)}}\gamma_{ij}(x))$, where $\sigma :TM\backslash\{
\hbox{Hyperplane}\}\to R$ is a smooth function. Using the ideas exposed in
\cite{2}, in these spaces, we can construct a Lagrange geometry and field
theory. This geometrical Lagrange theory will be regarded as a natural one
associated to the PDE system $\delta f=T$, in the sense of the first Udri\c ste's
question.

Now, we assume that a generalized Lagrange space $(M^n,g_{ij}(x,y))$
satisfies the following axioms:

a. 1. The fundamental tensor field $g_{ij}(x,y)$ is of the form
$$
g_{ij}(x,y)=e^{2\sigma(x,y)}\gamma_{ij}(x).
$$

a. 2. The space is endowed with the non-linear connection
$$
N^i_j(x,y)=\Gamma^i_{jk}(x)y^k,
$$
where $\Gamma^i_{jk}(x)$ are the Christoffel symbols for the Riemannian
metric $\gamma_{ij}(x)$.

Under these assumptions, our space verifies a constructive-axiomatic formulation of General Relativity
due to Ehlers, Pirani and Schild \cite{2}. This space represents a convenient
relativistic model, since it has the same conformal and projective properties as the
Riemannian space $(M,\gamma_{ij})$.

In the Lagrangian theory of electromagnetism, the electromagnetic tensors
$F_{ij}$ and $f_{ij}$ are
$$\displaystyle{
F_{ij}=\left(g_{ip}{\delta\sigma\over\delta x^j}-g_{jp}{\delta\sigma\over
\delta x^i}\right)y^p\;,\;f_{ij}=\left(g_{ip}{\partial\sigma\over\partial y^j}-
g_{jp}{\partial\sigma\over\partial y^i}\right)y^p
}.$$

Developping the formalism presented in \cite{2}, \cite{3} and denoting by
$r^i_{jkl}$ the curvature tensor field of the metric $\gamma_{ij}(x)$ , the
following Maxwell equations of the electromagnetic tensors hold
$$\left\{\begin{array}{lll}
\displaystyle{F_{ij\vert k}+F_{jk\vert i}+F_{ki\vert j}=-\sum_{(ijk)}g_{ip}
r^h_{qjk}{\partial\sigma\over\partial y^h}y^py^q},\\
F_{ij}\vert_k+F_{jk}\vert_i+F_{ki}\vert_j=-(f_{ij\vert k}+
f_{jk\vert i}+f_{ki\vert j}),\\
f_{ij}\vert_k+f_{jk}\vert_i+f_{ki}\vert_j=0,
\end{array}\right.
$$
where ${_\vert}_i$ (resp. $\vert_a$) represents the $h-$ (resp. $v-$) covariant
derivative induced by the non-linear connection $N^i_j$.

Suppose $\sigma=\sigma(x)$. Then the $v-$ electromagnetic tensor is
$f_{ij}=0$, the $h-$ covariant operator "${_\vert}_i$" becomes the
covariant derivative with respect to Levi-Civita connection of the metric
$g_{ij}(x,y)=e^{2\sigma(x)}\varphi_{ij}(x)$, the $h-$ electromagnetic tensor
$F_{ij}$ is the same with the classical electromagnetic tensor and the
Maxwell's equations reduce to the classical ones.

In the construction of the gravitational field equations, we shall use the
notations
$$\left\{\begin{array}{l}\medskip
r_{ij}=r^k_{ijk},\;r=\gamma^{ij}r_{ij},\;\displaystyle{
{\delta\over\delta x^i}={\partial\over \partial x^i}-
N^j_i{\partial\over\partial y^j}},\\
\sigma^{\scriptscriptstyle {H}}=
\displaystyle{\gamma^{kl}{\delta\sigma\over\delta x^k}{\delta\sigma\over\delta
x^l},\;\sigma^{\scriptscriptstyle {V}}=\gamma^{ab}{\partial\sigma\over\partial
y^a}{\partial\sigma\over\partial y^b},\;\overline\sigma=\gamma^{ij}\sigma_{ij},
\;\dot\sigma=\gamma^{ab}\dot\sigma_{ab}},
\end{array}\right.
$$
where
$\displaystyle{
\sigma_{ij}={\delta\sigma\over\delta x^i}\vert_j+{\delta\sigma\over\delta x^i}{\delta\sigma\over
\delta x^j}-{1\over 2}\gamma_{ij}\sigma^{\scriptscriptstyle H}},\;
\displaystyle{
\dot\sigma_{ab}=\left.{\partial\sigma\over\partial y^a}\right\vert_b+
{\partial\sigma\over\partial y^a}{\partial\sigma\over\partial y^b}-{1\over 2}
\gamma_{ab}\sigma^{\scriptscriptstyle V}}.$\\
In these conditions, the Einstein's equations of the space $(M,g_{ij}(x,y))$ take the form
$$\left\{
\begin{array}{ll}\medskip
\displaystyle{r_{ij}-{1\over 2}r\gamma_{ij}+t_{ij}={\cal K}
T^{\scriptscriptstyle H}_{ij}}\\
(2-n)(\dot\sigma_{ab}-\dot\sigma\gamma_{ab})={\cal K}T^{\scriptscriptstyle V}
_{ab},\\
\end{array}\right.
$$
where $T^{\scriptscriptstyle H}_{ij}$ and $T^{\scriptscriptstyle V}_{ab}$ are
the $h-$ and the $v-$ components of the energy momentum tensor field, ${\cal K}$
is the gravific constant and
$$
t_{ij}=(n-2)(\gamma_{ij}\overline\sigma-\sigma_{ij})+\gamma_{ij}r_{st}y^s
\gamma^{tp}{\partial\sigma\over\partial y^p}+{\partial\sigma\over\partial y^i}
r^a_{tja}y^t-\gamma_{is}\gamma^{ap}{\partial\sigma\over\partial y^p}
r^s_{tja}y^t.
$$

It is clear that the general metric $g_{ij}(x,y)=e^{2\sigma(x,y)}\varphi
_{ij}(x)$ implies Einstein equations which differ from the classical ones by the
additional tensor $t_{ij}$.

Finally, we remark that, in certain particular cases,
it is posible to build a generalized Lagrange geometry and field theory
naturally attached to a system of partial differential equations.
In these geometrical structures, the solutions of $C^2$ class of the PDE
system become harmonic maps. This idea was suggested by Udri\c ste in private discussions and in \cite{9},
\cite{10}.

{\bf Open problem.} Because the generalized Lagrange structure constructed
in this paper is not unique, it arises a natural question:

$-$Is it possible to build a unique generalized Lagrange geometry naturally
asociated to a given PDEs system?

An answer to this question will be offered by author in a subsequent paper,
using a more general Lagrange geometry, naturally attached to a
multidimensional Lagrangian defined on the jet fibration of order one.

{\bf Acknowledgements.} I would like to express my gratitude to the reviewer
of Southeast Asian Bulletin of Mathematics and Prof. Dr. C. Udri\c ste for
their valuable comments and very useful suggestions.

\begin{center}
University POLITEHNICA of Bucharest\\
Department of Mathematics I\\
Splaiul Independentei 313\\
77206 Bucharest, Romania\\
e-mail:mircea@mathem.pub.ro
\end{center}

\end{document}